\documentclass[10pt,twoside]{article}
\usepackage{graphicx}
\usepackage{amsmath}
\usepackage{amsfonts}
\usepackage{amssymb}
\usepackage{Latex-document}

\newcommand{\f}{\frac}
\newcommand{\Z}{\mathbb{Z}}
\newcommand{\R}{\mathbb{R}}
\newcommand{\on}{\operatorname}
\newcommand{\ca}{\mathcal}
\newcommand{\mf}{\mathfrak}
\newcommand{\Cl}{\on{Cl}}
\renewcommand{\o}{\mathfrak{o}}
\renewcommand{\S}{\ca{S}}
\newcommand{\ad}{\on{ad}}
\newcommand{\g}{\mf{g}}
\renewcommand{\d}{\on{d}}
\newcommand{\p}{\partial}

\newcommand{\W}{\ca{W}}
\newcommand{\wh}{\widehat}
\newcommand{\A}{\ca{A}}
\newcommand{\Om}{\Omega}

\title{\bf Clifford Algebras  and the \vskip -2mm
Duflo Isomorphism \vskip 6mm}
\author{E. Meinrenken\vspace*{-0.5cm}\thanks{Department of Mathematics,
University of Toronto, 100 St. George Street, Toronto, ON M6R1G7, Canada. E-mail: mein@math.toronto.edu}}
\date{\vspace{-8mm}}

\begin{document}

\maketitle

\thispagestyle{first} \setcounter{page}{637}

\begin{abstract}

\vskip 3mm

This article summarizes joint work with A. Alekseev (Geneva) on
the Duflo isomorphism for quadratic Lie algebras. We describe a
certain quantization map for Weil algebras, generalizing both the
Duflo map and the quantization map for Clifford algebras. In this
context, Duflo's theorem generalizes to a statement in equivariant
cohomology.

\vskip 4.5mm


\end{abstract}

\vskip 12mm

\section{Introduction} \label{section 0}\setzero

\vskip-5mm \hspace{5mm}

The universal enveloping algebra $U(\g)$ of a Lie algebra
$(\g,[\cdot,\cdot]_\g)$ is the quotient of the tensor algebra $T(\g)$
by the relations, $\xi\xi'-\xi'\xi=[\xi,\xi']_\g$. The inclusion of
the symmetric algebra $S(\g)$ into $T(\g)$ as totally symmetric
tensors, followed by the quotient map, gives an isomorphism of
$\g$-modules
\begin{equation}\label{eq:sym}
\on{sym}:\,S(\g)\to U(\g)
\end{equation}
called the {\em symmetrization map}. The restriction of $\on{sym}$ to
$\g$-invariants is a vector space isomorphism, but not an algebra
isomorphism, from invariant polynomials to the center of the
enveloping algebra.  Let $J\in C^\infty(\g)$ be the function
$$ J(\xi)=\det(j(\ad_\xi)),\ \ j(z)=\f{\on{sinh}(z/2)}{z/2},$$
and $J^{1/2}$ its square root (defined in a
neighborhood of $\xi=0$). Denote by
$\widehat{J^{1/2}}$ the infinite order differential operator on
$S\g\subset C^\infty(\g^*)$, obtained by replacing the variable $\xi\in\g$
with a directional derivative $\f{\p}{\p \mu}$, where $\mu$ is the
dual variable on $\g^*$. Duflo's celebrated theorem says that the composition
$$\on{sym}\circ  \widehat{J^{1/2}}:\ S\g\to U(\g)$$
restricts to an algebra isomorphism, $(S\g)^\g\to \on{Cent}(U(\g))$.
In more geometric language, Duflo's theorem gives an isomorphism
between the algebra of invariant constant coefficient differential
operators on $\g$ and bi-invariant differential operators on the
corresponding Lie group $G$.

The purpose of this note is to give a quick overview of
joint work with A. Alekseev \cite{am,am1}, in which
we obtained a new proof and a
generalization of Duflo's theorem for the special case of a {\em
quadratic} Lie algebra. That is, we assume that
$\g$ comes equipped with an invariant, non-degenerate, symmetric bilinear
form $B$. Examples of quadratic Lie algebras
include semi-simple Lie algebras, or the semi-direct
product $\g=\mf{s}\ltimes\mf{s}^*$ of a Lie algebra $\mf{s}$ with its
dual. Using $B$ we can define the Clifford algebra $\Cl(\g)$. Duflo's
factor $J^{1/2}(\xi)$
arises as the Berezin integral of $\exp(q(\lambda(\xi)))\in \Cl(\g)$,
where $q:\,\wedge(\g)\to \Cl(\g)$ is the quantization map, and
$\lambda:\,\g\to \wedge^2\g$ is the map dual to the Lie bracket.

\section{Clifford algebras} \label{section 1}\setzero
\vskip-5mm \hspace{5mm }

Let $V$ be a finite-dimensional real vector space, equipped with
a non-degenerate symmetric bilinear form $B$. Fix a basis $e_a\in V$
and let $e^a\in V$ be the dual basis. We denote by
$\o(V)\subset \on{End}(V)$ the space of endomorphisms $A$ of $V$ that
are skew-symmetric with respect to $B$. For any $A\in \o(V)$ we denote
its components by $A_{ab}=B(e_a,A e_b)$.
Consider the function $\ca{S}:\,\o(V)\to \wedge(V)$ given by
$$ \ca{S}(A)={\det}^{1/2}\left(j(A)\right)
\ \exp_{\wedge(V)}\left (\f{1}{2} f(A)_{ab}e^a\wedge e^b\right)$$
(using summation convention), where
\begin{equation}\label{eq:log}
f(z)=(\ln j)'(z)=\frac{1}{2} \coth(\f{z}{2}) - \f{1}{z}.
\end{equation}
In turns out that, despite the
singularities of the exponential, $\ca{S}$ is a global analytic
function on all of $\o(V)$. It has the following nice property. Let
$\Cl(V)$ denote the Clifford algebra of $V$, defined as a quotient of
the tensor algebra $T(V)$ by the relations $vv'+v'v=B(v,v')$. The
inclusion of $\wedge(V)$ into $T(V)$ as totally anti-symmetric tensors,
followed by the quotient map to $\Cl(V)$, gives a vector space
isomorphism
$$ q:\,\,\wedge(V)\to \Cl(V)$$
known as the {\em quantization map}. Then $\S(A)$ relates the exponentials
of quadratic elements $1/2 A_{ab}e^a \wedge e^b$
in the exterior algebra with the exponentials
of the corresponding elements  $1/2 A_{ab}e^a e^b$
in the Clifford algebra:
\begin{equation}\label{eq:expo}
 \exp_{\Cl(V)}(1/2 A_{ab}e^a e^b)=q\left(\iota(\S(A))
\exp_{\wedge(V)}(1/2 A_{ab}e^a\wedge  e^b)\right).
\end{equation}
Here $\iota:\,\wedge(V)\to \on{End}(V)$ is the contraction operator.
In fact, one may add linear terms to the exponent: Let $E$ be some
vector space of ``parameters'', and $\phi^a\in E$. Then the following identity
holds in the $\Z_2$-graded tensor product $\Cl(V)\otimes \wedge(E)$:
$$ \exp(1/2 A_{ab}e^a e^b+e_a\otimes \phi^a)=q\Big(\iota(\S(A))
\exp(1/2 A_{ab}e^a \wedge e^b+e_a\otimes \phi^a)\Big).
$$

\section{Quadratic Lie algebras} \label{section 2}
\setzero\vskip-5mm \hspace{5mm }

Let us now consider the case $V=\g$ of a quadratic Lie algebra.
Invariance of the bilinear form $B$  means that the the adjoint representation
$\ad:\,\g\to \on{End}(\g)$ takes values in $\o(\g)$, or equivalently
that the structure constants $f_{abc}=B(e_a,[e_b,e_c])$ are invariant
under cyclic permutations of the indices $a,b,c$.  We specialize
\eqref{eq:expo} to  $A=\ad_\xi$ for $\xi\in\g$, so that
$\lambda^\g:\,\g\to \wedge^2\g$,
$$ \lambda^\g(\xi)=1/2 (\ad_\xi)_{ab}e^a\wedge e^b$$
is the map dual to the Lie bracket. Also, take $E=T_\xi^*\g$
and $\phi^a=-d\, \xi^a$, where $\xi^a=B(\xi,e^a)$
are the coordinate functions. Then our formula become the following
identity in $\Cl(\g)\otimes \Omega(\g)$
\begin{equation}\label{eq:2}
\exp\big(q(\lambda^\g)-e_a\, d\xi^a\big)
=q\Big(\iota(\S^\g)\exp(\lambda^\g-e_a\, d\xi^a)\Big),
\end{equation}
where $\S^\g=\S\circ\ad:\,\g\to \wedge\g$. Consider now the following
cubic element in the Clifford algebra,
$$ \ca{C}=\f{1}{6}f_{abc}e^a\, e^b \, e^c\in\Cl(\g).$$
A beautiful observation of Kostant-Sternberg
\cite{ks} says that $\ca{C}$ squares to a {\em constant}:
$$ \ca{C}^2=-\f{1}{48} f_{abc} f^{abc}.$$
It follows that the graded commutator $\d^{\on{Cl}\g}:=
[\ca{C},\cdot\,]$ defines a differential on $\Cl(\g)$.
This {\it Clifford differential}
is compatible with the filtration of $\Cl(\g)$, and the induced
differential $\d^{\wedge \g}$
on the associated graded algebra
$\on{gr}(\Cl(\g))=\wedge\g$ is nothing but the Lie algebra
differential. Let $\d^{\on{Rh}}$ denote the exterior differential
on the deRham complex $\Omega(\g)$.

It is easily verified
that $\lambda^\g-e_a\, d\xi^a\in\wedge \g\otimes \Omega(\g)$ is
closed for the differential $\d^{\wedge \g}+\d^{\on{Rh}}$,
while $q(\lambda^\g)-e_a\, d\xi^a\in \Cl(\g)\otimes \Omega(\g)$
is closed under $\d^{\Cl(\g)}+\d^{\on{Rh}}$. Together with
\eqref{eq:2}, this leads to a number of consistency conditions for
the function $\S^\g$. One of these conditions gives a solution
of the {\em classical dynamical Yang-Baxter equation} (CDYBE):  Let
$\mf{r}:\,\g\to\o(\g)$ be the meromorphic function
$ \mf{r}^\g(\xi)=f(\ad_\xi)$ appearing in the exponential
 factor of $\S^\g$. Then
\begin{equation}
 \on{cycl}_{abc}\left(\f{\p \mf{r}_{ab}}{\p \xi^c}-
\mf{r}_{ak}f^{kl}_b\mf{r}_{lc}\right)=-\f{1}{4}f_{abc}
\end{equation}
where $\on{cycl}_{abc}$ denotes the sum over cyclic permutations of
$a,b,c$. This solution of the CDYBE was obtained by Etingof-Varchenko
\cite{ev} and in \cite{am} by different methods. In Etingof-Schiffmann
\cite{es}, it is shown that $\mf{r}_{ab}$ is in fact the {\em unique}
solution of this particular CDYBE, up to gauge transformation.
More general CDYBE's are associated to a
pair $\mf{h}\subset \g$ of Lie algebras, here $\mf{h}=\g$. The proof sketched
above can be modified to produce some of these more general solutions.

\section{The non-commutative Weil algebra}
\label{section 3} \setzero\vskip-5mm \hspace{5mm }

Using $B$ to identify the Lie algebra $\g$ with its dual $\g^*$, the
Weil algebra of $\g$ is the $\Z$-graded $\g$-module given as a
tensor product
$$ W\g=S\g\otimes \wedge\g,$$
where generators of $S\g$ are assigned degree $2$. Let $L_\xi^W$
for $\xi\in\g$
denote the generators for the $\g$-action on $W\g$, and
$\iota_\xi^W=1\otimes\iota_\xi$ the contraction operators.  The {\it Weil
differential} $\d^W$ is a derivation of degree 1,
uniquely characterized by its properties
$\d^W\circ\d^W=0$ and
$\d^W(1\otimes\xi)=\xi\otimes 1$ for $\xi\in \g$.
The Weil algebra $W\g$ with these three types of derivations is an
example of a $\g$-differential algebra: That is,
$L_\xi^W,\iota_\xi^W,\d^W$ satisfy relations similar to contraction
operators, Lie derivatives, and de Rham differential for a manifold
with group action.

In \cite{am}, we introduced the following non-commutative version of
the Weil algebra,
$$ \ca{W}\g=U\g\otimes \Cl(\g).$$
It carries a $\Z$-filtration, where generators of $U(\g)$ are assigned
filtration degree 2, with associated graded algebra $\on{gr}(\W\g)=W\g$.
Moreover, it carries a $\Z_2$-grading, compatible with the $\Z$-filtration
in the sense of \cite{ks}.
Define contraction operators as $\Z_2$-graded commutators
$\iota_\xi^\W=[1\otimes\xi,\cdot\,]$, let
$L_\xi^\W$ be the generators for the natural $\g$-module structure,
and set $\d^\W=[\ca{D},\cdot]$ where
$$ \ca{D}=e_a\otimes e^a-1\otimes \ca{C}\in \W\g$$
is the {\em cubic Dirac operator} \cite{ko}. Its square
$$ \ca{D}^2=\f{1}{2} e_a e^a\otimes 1-\f{1}{48}f_{abc} f^{abc} $$
is in the center of $\W\g$, hence $\d^\W$ is a differential.
As it turns out, $\W\g$ is again a $\g$-differential algebra.
The derivations $\d^\W,\iota_\xi^\W,L_\xi^\W$ respect the
$\Z$-filtration, and the induced derivations on the associated
graded algebra are just the standard derivations for
the Weil algebra $W\g$.

The vector space isomorphism $\on{sym}\otimes q:\,W\g\to \W\g$
intertwines the contraction operators and Lie derivatives,
but not the differentials. There does exist, however, a better
``quantization map'' $\ca{Q}:\,W\g\to \W\g$ that is also a
chain map. Using our function
$\S^\g\in C^\infty(\g)\otimes\wedge\g$, let $\iota(\wh{\S^\g})$
denote the operator on $W\g$, where the $\wedge\g$-factor acts
by contraction on $\wedge\g$ and the $C^\infty(\g)$-factor
as an infinite order differential operator.
\vskip.3in

\noindent {\bf Theorem.} \cite{am} {\it The quantization map
$$ \ca{Q}:=(\on{sym}\otimes q)\circ \iota(\wh{\S^\g}):\,W\g\to \W\g$$
intertwines the contraction operators, Lie derivatives, and
differentials on $W\g$ and on $\W\g$.}

The fact that $\ca{Q}$ intertwines the two differentials
$\d^W,\ \d^\W$ relies on  a number of special properties of the
function $\S^\g$, including the CDYBE.

Put differently, the quantization map $ \ca{Q}$ defines a new, graded
non-commutative ring structure on the Weil algebra $W\g$, in
such a way that the derivations $\iota_\xi^W,L_\xi^W,\d^W$
are still derivations for the new ring structure, and in fact become
{\em inner} derivations. Notice that $\ca{Q}$ restricts to the
quantization map for Clifford algebras
$q:\,\wedge\g\to \Cl(\g)$ on the second factor and to the
Duflo map on the first factor, but is not just the product
of these two maps.


\section{Equivariant cohomology}
\label{section 4} \setzero\vskip-5mm \hspace{5mm }

H. Cartan in \cite{ca} introduced the Weil algebra $W\g$ as an algebraic model
for the algebra of differential forms on the classifying bundle $EG$,
at least in the case $G$ compact.

In particular, it can be used to compute the
equivariant cohomology $H_G(M)$ (with real coefficients) for any
$G$-manifold $M$. Let $\iota_\xi^{Rh},L_\xi^{Rh},
\d^{Rh}$ denote the contraction operators, Lie derivatives, and differential
on the de Rham complex $\Omega(M)$ of differential forms. Let
$$ H_\g(M)=H((W\g\otimes \Om(M))_{\on{basic}},\d^W+\d^{Rh})$$
where $(W\g\otimes \Om(M))_{\on{basic}}$ is the subspace
annihilated by all Lie derivatives $L_\xi^W+L_\xi^{Rh}$ and all
contraction operators $\iota_\xi^W+\iota_\xi^{Rh}$. Cartan's
result says that $H_\g(M)=H_G(M,\R)$ provided $G$ is compact.

More generally, we can define $H_\g(\A)$ for any $\g$-differential
algebra $\A$. Let $\ca{H}_\g(\A)$ be defined by replacing $W\g$
with $\W\g$. The quantization map $\ca{Q}:\,W\g\to \W\g$ induces
a map $\ca{Q}:\,H_\g(\A)\to \ca{H}_\g(\A)$.

\noindent {\bf Theorem.}  \cite{am} {\it For any $\g$-differential algebra $\ca{A}$, the vector space isomorphism
$\ca{Q}:\,H_\g(\A)\to \ca{H}_\g(\A)$ is in fact an algebra isomorphism.}

Our proof is by construction of an explicit chain homotopy
between the two maps $W\g\otimes W\g\to \W\g$ given by
``quantization followed by multiplication'' and ``multiplication
followed by quantization'', respectively. Taking $\A$ to be the
trivial $\g$-differential algebra (i.e. $\A=\Om(\on{point})$),
the statement specializes to Duflo's theorem for quadratic $\g$.

\label{lastpage}

\end{document}